\documentclass{amsart2000}
\usepackage{hyperref}

\theoremstyle{plain}
\newtheorem{thm}{Theorem}
\newtheorem{prop}[thm]{Proposition}
\newtheorem{lem}[thm]{Lemma}
\newtheorem{cor}[thm]{Corollary}

\theoremstyle{remark}
\newtheorem{rem}[thm]{Remark}
\newtheorem{exm}[thm]{Example}

\renewcommand{\b}[1]{\mathbf{#1}}

\newcommand{\A}{{\mathcal A}}
\newcommand{\LL}{{\mathcal L}}
\newcommand{\B}{{\mathcal B}}
\newcommand{\F}{{\mathcal F}}
\newcommand{\HH}{{\mathcal H}}

\newcommand{\GG}{{\mathcal G}}

\newcommand{\bH}{{\mathbb H}}
\newcommand{\C}{{\mathbb C}}
\newcommand{\bP}{{\mathbb P}}
\newcommand{\Z}{{\mathbb Z}}
\newcommand{\II}{{\mathbb I}}

\newcommand{\sfD}{{\mathsf D}}

\newcommand\la{{\lambda}}

\DeclareMathOperator{\GL}{GL}
\DeclareMathOperator{\End}{End}
\DeclareMathOperator{\ii}{i}
\DeclareMathOperator{\Perv}{Perv}
\DeclareMathOperator{\codim}{codim}

\begin{document}

\title[Nonresonance conditions for arrangements]
{Nonresonance conditions for arrangements}
\author[D.~Cohen]{Daniel C.~Cohen$^\dag$}
\address{Department of Mathematics, Louisiana State University,
Baton Rouge, LA 70803}
\email{\href{mailto:cohen@math.lsu.edu}{cohen@math.lsu.edu}}
\urladdr{\href{http://www.math.lsu.edu/~cohen/}
{http://www.math.lsu.edu/\~{}cohen}}
\thanks{{$^\dag$}Partially supported by Louisiana Board of Regents
grant LEQSF(1999-2002)-RD-A-01 and by National Security Agency grant
MDA904-00-1-0038}

\author[A. Dimca]{Alexandru Dimca}
\address{Laboratoire de Math. Pures, Universit\'e Bordeaux I, France}
\email{\href{mailto:dimca@math.u-bordeaux.fr}{dimca@math.u-bordeaux.fr}}
\urladdr{\href{http://www.math.u-bordeaux.fr/~dimca/}
{http://www.math.u-bordeaux.fr/\~{}dimca}}

\author[P.~Orlik]{Peter Orlik$^\ddag$}
\address{Department of Mathematics, University of Wisconsin,
Madison, WI 53706}
\email{\href{mailto:orlik@math.wisc.edu}{orlik@math.wisc.edu}}
\thanks{{$^\ddag$}Partially supported by National Security Agency
grant MDA904-02-1-0019}

% 32S22 Relations with arrangements of hyperplanes
%       (Several complex variables and analytic spaces; Singularities)
% 52C35 Arrangements of points, flats, hyperplanes
%       (Convex and discrete geometry; Discrete geometry)
% 55N25 Homology with local coefficients, equivariant cohomology
%       (Algebraic topology; Homology and cohomology theories)


\begin{abstract}
We prove a vanishing theorem for the cohomology of the complement of a
complex hyperplane arrangement with coefficients in a complex local
system.  This result is compared with other vanishing theorems, and
used to study Milnor fibers of line arrangements, and hypersurface
arrangements.
\end{abstract}

%\date{October 26, 2002}

\maketitle

\section{Introduction}
\label{sec:intro}

Let $\A$ be an arrangement of hyperplanes in the complex projective
space $\bP^n$, with complement $M(\A)=\bP^n \setminus
\bigcup_{H\in\A}H$.  Let $\LL$ be a complex local system of
coefficients on $M(\A)$.  The need to calculate the local system
cohomology $H^*(M(\A),\LL)$ arises in a variety of contexts, including
the Aomoto-Gelfand theory of multivariable hypergeometric integrals
\cite{AK,Gel1}; representation theory of Lie algebras and quantum
groups and solutions of the Knizhnik-Zamolodchikov differential
equation in conformal field theory \cite{Va}; and the determination of
the cohomology groups of the Milnor fiber of the non-isolated
hypersurface singularity at the origin in $\C^{n+1}$ associated to the
arrangement $\A$ \cite{CS}.

In light of these applications, and others, the cohomology
$H^*(M(\A),\LL)$ has been the subject of considerable recent interest. 
Call the local system $\LL$ {\em nonresonant} if this cohomology is
concentrated in dimension $n$, that is, $H^k(M(\A),\LL)=0$ for $k \neq
n$.  Necessary conditions for vanishing, or nonresonance, have been
determined by a number of authors, including Esnault, Schectman, and
Viehweg \cite{ESV}, Kohno~\cite{Ko}, and Schechtman, Terao, and
Varchenko \cite{STV}.  Many of these results make use of Deligne's
work \cite{De}, and thus require the realization of $M(\A)$ as the
complement of a normal crossing divisor in a complex projective
manifold.

An edge is a nonempty intersection of hyperplanes.  An edge is {\em
dense} if the subarrangement of hyperplanes containing it is
irreducible: the hyperplanes cannot be partitioned into nonempty sets
so that after a change of coordinates hyperplanes in different sets
are in different coordinates.  This is a combinatorially determined
condition which can be checked in a neighborhood of a given edge, see
\cite{STV}.  Consequently, this notion makes sense for both affine and
projective arrangements.  Let $\sfD(\A)$ denote the set of dense edges
of the arrangement $\A$.

Let $N=\bigcup_{H\in\A} H$ be the union of the hyperplanes of $\A$. 
There is a canonical way to obtain an embedded resolution of the
divisor $N$ in $\bP^n$.  First, blow up the dense 0-dimensional edges
of $\A$ to obtain a map $p_1:Z_1 \to \bP^n$.  Then, blow up all the
proper transforms under $p_1$ of projective lines corresponding to
dense 1-dimensional edges in $\sfD(\A)$.  Continuing in this way, we
get a map $p=p_{n-1}:Z_{n-1} \to \bP^n$ which is an embedded
resolution of the divisor $N$ in $\bP^n$.  Let $Z=Z_{n-1}$.  Then,
$D=p^{-1 }(N)$ is a normal crossing divisor in $Z$, with smooth
irreducible components $D_X$ corresponding to the edges $X \in
\sfD(\A)$.  Furthermore, the map $p$ induces a diffeomorphism $Z
\setminus D =M(\A)$, see \cite{OT2,STV,Va} for details.

Let $\LL$ be a complex local system of rank $r$ on the complement
$M(\A)$ associated to a representation
\[
\rho : \pi _1 (M(\A ),a) \to \GL_r(\C).
\]
To each irreducible component $D_X$ of the normal crossing divisor $D$
corresponds a well-defined conjugacy class $T_X$ in $\GL_r(\C)$,
obtained as the monodromy of the local system $\LL$ along a small loop
turning once in the positive direction about the hypersurface $D_X$. 
In this note, we prove the following nonresonance theorem.

\begin{thm} \label{thm:nonres}
Assume that there is a hyperplane $H \in \A$ such that for any dense
edge $X \in \sfD(\A)$ with $X \subseteq H$ the corresponding monodromy
operator $T_X$ does not admit $1$ as an eigenvalue.  Then 
$H^k(M(\A ),\LL)=0$ for any $k \ne n$.
\end{thm}

In the case when $\LL$ is a rank one complex local system arising in
the context of the Milnor fiber associated to $\A$ (see \cite{CS} and
Section \ref{sec:mf} below), this result was obtained by Libgober
\cite{Li}.  We do not see a simple way to extend the topological proof
given by Libgober in this special case to the general case stated
above.

The structure of this note is as follows.  Theorem \ref{thm:nonres} is
proved in Section \ref{sec:proof}.  Local systems which arise from
flat connections on trivial vector bundles are considered in Section
\ref{sec:flat}.  The implications of Theorem \ref{thm:nonres} in this
special case are compared with other nonresonance theorems in Section
\ref{sec:compare}.  A brief application to Milnor fibers associated to
line arrangements in $\bP^2$, strengthening a result of Massey
\cite{Ma}, is given in Section~\ref{sec:mf}.  A strategy to handle
arrangements of more general hypersurfaces is presented in
Section~\ref{sec:gen}.

\section{Proof of Theorem \ref{thm:nonres}}
\label{sec:proof}

Before giving the proof of Theorem \ref{thm:nonres}, we say a few
words concerning the monodromy operators $T_X$.  Let $\widetilde\A$
denote the central arrangement in $\C^{n+1}$ corresponding to the
arrangement $\A$ in $\bP^n$.  The representation $\rho:\pi_1(M(\A),a)
\to\GL_r(\C)$ induces a representation
\[
\tilde\rho: \pi _1 (M(\widetilde\A),\tilde{a}) \to GL_r(\C),
\]
where $\tilde{a}$ is any lift of the base point $a$.

For such a representation $\tilde\rho$, associated to a local system
$\widetilde\LL$ on the complement of $\widetilde\A$, there is a
well-defined {\em total turn monodromy operator}
$T(\widetilde\A)=\tilde\rho(\gamma)$, where
$\gamma(t)=\exp(2 \pi \ii t)\tilde{a}$ for $t \in [0,1]$.  Choosing a
generic line passing through $\tilde{a}$ and close to (but not
through) the origin yields, in the usual way, $m=|\widetilde\A|$
elementary loops which generate the fundamental group
$\pi_1(M(\widetilde\A),\tilde{a})$.  The product in a certain natural
order of the associated monodromy operators in $\GL_r(\C)$ is easily
seen to be exactly the total turn monodromy operator
$T(\widetilde\A)$.  In particular, if $\widetilde\LL$ (resp.,
$\tilde\rho$) is abelian, the order is irrelevant and
$T(\widetilde\A)=T_1\cdots T_m$, where $T_j$ is the monodromy
around the hyperplane $H_j$.

Let $X$ be a dense edge of $\widetilde{\A}$ and let $V$ be an affine
subspace of complementary dimension in $\C^{n+1}$, transverse to $X$
at a generic point $x$ in $X$.  Consider $V$ as a vector space with
origin $x$ and let $\widetilde\B$ denote the central arrangement in
$V$ induced by $\widetilde\A$.  Let $\B$ be the corresponding
projective arrangement in $\bP(V)$.  The following result can be
checked by the reader.

\begin{lem} \label{lem:cong}
The following conjugacy classes in $\GL_r(\C)$ coincide
\begin{enumerate}
\item[(a)] The monodromy $T_X$;
\item[(b)] The total turn monodromy $T(\widetilde\B)$;
\item[(c)] The monodromy of the restriction of the local system $\LL$
to any of the fibers of the projection $M(\widetilde\B) \to M(\B)$.
\end{enumerate}
\end{lem}

Now we prove Theorem \ref{thm:nonres}.  In this proof, we use a
partial resolution similar to the resolution $p: Z \to \bP^n$
mentioned in the Introduction.

First, blow up all the dense 0-dimensional edges contained in the
chosen hyperplane $H$.  This yields a map $q_1:W_1 \to \bP^n$.  Then,
blow up all the proper transforms under $q_1$ of projective lines
corresponding to dense 1-dimensional edges in $\sfD(\A)$ which are
contained in $H$.  Continuing in this way, we get an embedded
resolution of the divisor $N$ in $\bP^n$ along $H$.  This yields a map
$q=q_{n-1}:W=W_{n-1} \to \bP^n$ such that $E=q^{-1}(N)$ is a normal
crossing divisor at any point of $H'=q^{-1}(H)$.  Moreover, $H'$ has
smooth irreducible components $E_X$ corresponding to the edges $X \in
\sfD(\A)$, $X \subseteq H$, and $q$ induces a diffeomorphism $W
\setminus H' =\bP^n \setminus H$.  Note that the conjugacy classes
$T_X$ for $X \in \sfD(\A)$, $X \subseteq H$ constructed from the
resolutions $Z$ and $W$ coincide.

Let $U=W \setminus H'= \bP^n \setminus H$, and let $i: M(\A) \to U$
and $j:U \to W$ be the corresponding inclusions.  Denote the derived
category of constructible bounded complexes of sheaves of $\C$-vector
spaces on $U$ by $\b{D}_c^b(U)$, and let $\F=Ri_* \LL[n] \in
\b{D}_c^b(U)$.  Since $M(\A)$ is smooth, the shifted local system
$\LL[n]$ is a perverse sheaf on $ M(\A)$.  Moreover, since $i$ is a
quasi-finite affine morphism, it preserves perverse sheaves, and hence
$\F \in \Perv(U)$, see \cite[(10.3.27)]{KS}.

Since $U$ is a smooth affine variety, by the Artin Theorem and Verdier
Duality, see \cite[(10.3.5) and (10.3.8)]{KS}, we have the following
vanishing results:
\begin{equation} \label{eqn:vanishing}
\bH^k(U,\F)=0\text{ for all $k>0$, and } \bH_c^k(U,\F)=0
\text{ for all $k<0$.}
\end{equation}
Note that we can write $\bH^k(U,\F)=H^k(Rs_*Rj_*\F)$ and
$\bH_c^k(U,\F)=H^k(Rs_!Rj_!\F)$, where $s:W \to {\rm pt}$ is the
constant map to a point.  Since $\bP^n$ is compact, $s$ is a proper
map, and hence $Rs_*=Rs_!$.  Consequently, in light of the Leray-type
isomorphism
\[
H^{k+n}(M(\A),\LL)= \bH^k(U,\F),
\]
to prove Theorem \ref{thm:nonres}, it is enough to establish the
following result.

\begin{lem} \label{lem:canonical}
With the above notation, if for any dense edge $X \in \sfD(\A)$
with $X \subseteq H$ the corresponding monodromy operator $T_X$
does not admit 1 as an eigenvalue, then the canonical morphism
$Rj_!\F \to Rj_*\F$ in $\b{D}_c^b(\bP^n )$ is an isomorphism.
\end{lem}
\begin{proof}
The canonical morphism is an isomorphism if and only if the induced
morphisms on the level of stalk cohomology are isomorphisms.  This
local property is clearly satisfied for the stalks at $x \in U$ since
$U$ is open.

Consider the case $x \in H'$.  Then $\HH^*(Rj_!\F)_x=0$ using the
proper base change, see \cite[(2.6.7)]{KS}.  To show that
$\HH^*(Rj_*\F)_x=0$, we have to compute the cohomology groups
$\HH^k(Rj_*\F)_x=H^{k+n}(M(\A) \cap B, \LL)$, where $B$ is a small
open ball in $W$ centered at $x$.

Since $E$ is a normal crossing divisor at $x$, it follows that the
fundamental group of $ M(\A) \cap B= (W \setminus E) \cap B$ is
abelian.  Using the methods of \cite{EV}, we can decrease the rank of
the local system $\LL$.  Repeating this process yields a rank one
local system, where the result follows using the K\"unneth formula,
since at least one of the irreducible components of $E$ passing
through $x$ corresponds to a dense edge $X \subseteq H$.

This completes the proof of Lemma \ref{lem:canonical}, and hence that
of Theorem \ref{thm:nonres} as well.
\end{proof}

\begin{rem}
A vanishing result similar to Theorem \ref{thm:nonres} for hyperplane
arrangements over algebraically closed fields of positive
characteristic may be obtained by using \cite{BBD} as a reference
instead of \cite{KS}.
\end{rem}

\begin{rem}
Assume that there is a hyperplane $H \in \A$ such that for any dense
edge $X \in \sfD(\A)$ with $X \subseteq H$ and $\codim X \leq c$ the
corresponding monodromy operator $T_X$ does not admit 1 as an
eigenvalue.  Then $H^p(M(\A),\LL)=0$ for any $p$ with $0 \leq p < c$. 
Indeed, by intersecting with a generic affine subspace $E$ with $\dim
E=c$, we obtain a $c$-homotopy equivalence $M(\A) \cap E \to M(\A)$
induced by the inclusion, and hence isomorphisms $H^p(M(\A) \cap
E,\LL)=H^p(M(\A),\LL)$ for $0\leq p<c$.  The assertion follows by
applying Theorem \ref{thm:nonres} to the arrangement in $E$ induced by
the arrangement $\A$.
\end{rem}

\section{A Special Case}
\label{sec:flat}

In this section, we consider the special case of local systems which
arise from flat connections on trivial vector bundles.  Write
$\A=\{H_1,\dots,H_m\}$ and for each $j$, let $f_j$ be a linear form
with zero locus $H_j$.  Let $\omega_j=d\log(f_j)$, and choose $r\times
r$ matrices $P_j \in \End(\C^r)$ which satisfy $\sum_{j=1}^m P_j=0$. 
For an edge $X$ of $\A$, set $P_X=\sum_{X \subseteq H_j} P_j$. 
Consider the connection on the trivial vector bundle of rank $r$ over
$M(\A)$ with $1$-form $\omega=\sum_{j=1}^m \omega_j \otimes P_j$.  The
connection is flat if $\omega \wedge \omega =0$.  This is the case if
the endomorphisms $P_j$ satisfy
\begin{equation} \label{eqn:holonomy}
[P_j,P_X]=0\quad \text{for all $j$ and edges $X$ such that
$\codim X=2$ and $X \subseteq H_j$,}
\end{equation}
see \cite{Ko}.  Let $\LL$ be the rank $r$ complex local system on
$M(\A)$ corresponding to the flat connection on the trivial vector
bundle over $M(\A)$ with $1$-form $\omega$.

\begin{rem}\label{rem:compare}
An arbitrary local system $\LL$ on $M(\A)$ need not arise as the sheaf
of horizontal sections of a trivial vector bundle equipped with a flat
connection as described above.  The existence of such a connection is
related to the Riemann-Hilbert problem for $\LL$, see Beauville
\cite{Beau}, Bolibrukh \cite{Bo}, and Kostov \cite{Kos}.  Even in the
simplest case, when $n=1$ and $|\A|>3$, there are local systems $\LL$
of any rank $r \geq 3$ on $M(\A)$ for which the Riemann-Hilbert
problem has no solution, see \cite[Theorem~3]{Bo}.
\end{rem}

For a local system which may be realized as the sheaf of horizontal
sections of a trivial vector bundle equipped with a flat connection,
Theorem \ref{thm:nonres} has the following consequence.

\begin{cor} \label{cor:nonres1}
Assume that there is a hyperplane $H \in \A$ such that
\begin{equation} \label{eqn:nonres1}
\text{none of the eigenvalues of $P_X$ lies in $\Z$ for every
dense edge $X \subseteq H$.}
\end{equation}
Then
\[
H^k(M(\A),\LL)=0\text{ for }k\neq n.
\]
\end{cor}
This result is a refinement of the vanishing theorem of Kohno
\cite{Ko}, where condition \eqref{eqn:nonres1} is required to hold for
all edges.  Next, we recall the following well known nonresonance
theorem of Schechtman, Terao, and Varchenko \cite{STV}.

\begin{thm}[{\cite[Corollary 15]{STV}}] \label{thm:stv}
Assume that none of the eigenvalues of $P_X$ lies in $\Z_{\ge 0}$
for every dense edge $X \in \sfD(\A)$. Also suppose that
$P_iP_j=P_jP_i$ for all $i,j$.  Then
\[
H^k(M(\A),\LL)=0\text{ for }k\neq n.
\]
\end{thm}
Note that this result pertains only to abelian local systems. This
assumption is not necessary.

\begin{thm} \label{thm:stvImprovement}
Assume that
\begin{equation} \label{eqn:stvImproved}
\text{none of the eigenvalues of $P_X$ lies in $\Z_{\ge 0}$ for
every dense edge $X \in \sfD(\A)$.}
\end{equation}
Then
\[
H^k(M(\A),\LL)=0\text{ for }k\neq n.
\]
\end{thm}
\begin{proof}[Sketch of Proof]
Let $B^\bullet(\A)$ denote the algebra of global differential forms on
$M(\A)$ generated by the $1$-forms $\omega_j$, the Brieskorn algebra
of $\A$.  Since the endomorphisms $P_j$ satisfy \eqref{eqn:holonomy},
the tensor product $B^\bullet(\A)\otimes \C^r$, with differential
given by multiplication by $\omega=\sum_{j=1}^m \omega_j \otimes P_j$,
is a complex, which may be realized as a subcomplex of the twisted de
Rham complex of $M(\A)$ with coefficients in $\LL$.

By work of Esnault, Schectman, and Viehweg \cite{ESV}, refined by
Schechtman, Terao, and Varchenko \cite{STV}, the above assumptions on
the eigenvalues of $P_X$ imply that there is an isomorphism
\[
H^*(M(\A);\LL) \cong H^*(B^\bullet(\A),\omega).
\]
Thus it suffices to show that $H^q(B^\bullet(\A),\omega)=0$ for
$q\neq n$.

For an abelian local system, this was established by Yuzvinsky
\cite{Yuz}.  To extend his argument to an arbitrary local system, it
is enough to show that the complex $(B^\bullet(\A),\omega)$ is acyclic
for a central arrangement $\A$.  This may be accomplished using the
Euler derivation to produce a chain contraction, a straightforward
modification of the proof given by Yuzvinsky.
\end{proof}

\section{Comparison}
\label{sec:compare}

The purpose of this section is to compare the nonresonance results of
the previous section.  A local system (resp., a collection
$(P_1,\dots,P_m)$ of endomorphisms satisfying \eqref{eqn:holonomy} and
$\sum_{j=1}^m P_j=0$) will be called $\A$--{\em nonresonant} if it
satisfies condition \eqref{eqn:stvImproved}, and will be called
$(\A,H)$--{\em nonresonant} if it satisfies condition
\eqref{eqn:nonres1}.  Let $\II_r$ denote the $r\times r$ identity
matrix, and note that if $k_1,\dots,k_m$ are integers with
$\sum_{j=1}^m k_j=0$, then the collections of endomorphisms
$(P_1,\dots,P_m)$ and $(P_1+k_1\cdot \II_r, \dots,P_m+k_m\cdot \II_r)$
give rise to the same representation and local system.  Furthermore,
if the endomorphisms $P_j$ satisfy the conditions of
\eqref{eqn:holonomy}, then so do the endomorphisms
$P_j+k_j\cdot\II_r$.  Hence, if the connection with $1$-form
$\sum_{j=1}^m \omega_j \otimes P_j$ is flat, then so is the connection
with $1$-form $\sum_{j=1}^m \omega_j \otimes (P_j+k_j\cdot\II_r)$.

Surprisingly, the monodromy condition \eqref{eqn:nonres1} of
Corollary \ref{cor:nonres1} is more stringent than the condition
\eqref{eqn:stvImproved} of Theorem \ref{thm:stvImprovement}.

\begin{prop} \label{prop:Comparison}
Let $\LL$ be an rank $r$ complex local system on $M(\A)$ induced
by a collection of $(\A,H)$--nonresonant endomorphisms
$(P_1,\dots,P_m)$. Then there are integers $k_1,\dots,k_m$ so that
the collection $(P_1+k_1\cdot \II_r,\dots,P_m+k_m\cdot \II_r)$ is
$\A$--nonresonant.
\end{prop}
\begin{proof}
Without loss, assume that the collection of endomorphisms
$(P_1,\dots,P_m)$ is $(\A,H_1)$--nonresonant.  Then by
\eqref{eqn:nonres1}, for every dense edge $X$ of $\A$ for which $X
\subseteq H_1$, the eigenvalues of $P_X$ are not integers.

Let $q$ be a positive integer, greater than the maximum of the
absolute values of the eigenvalues of the endomorphisms $P_Y$, where
$Y$ ranges over all dense edges of $\A$ for which $Y \not\subseteq
H_1$.  Let
\[
\widehat{P}_j=
\begin{cases}
P_1+(m-1)\cdot q\cdot \II_r&\text{if $j=1$,}\\
P_j-q\cdot \II_r&\text{if $j \neq 1$,}
\end{cases}
\]
and note that $\sum_{j=1}^m \widehat{P}_j = \sum_{j=1}^m P_j = 0$.

We assert that $(\widehat{P}_1,\widehat{P}_2,\dots,\widehat{P}_m)$ is
an $\A$-nonresonant collection of endomorphisms.  For this, let $X$ be
a dense edge of $\A$.  If $X \subseteq H_1$, then the eigenvalues of
$\widehat{P}_X=P_X+(m-|X|)\cdot q\cdot\II_r$ are not integers since
the eigenvalues of $P_X$ are not integers.  If $X \not\subseteq H_1$,
then the eigenvalues of $\widehat{P}_X=P_X-|X|\cdot q \cdot\II_r$ are
not in $\Z_{\ge 0}$ by the choice of $q$.  Hence the collection of
endomorphisms $(\widehat{P}_1,\widehat{P}_2,\dots,\widehat{P}_m)$
satisfies \eqref{eqn:stvImproved}, and is thus $\A$--nonresonant.
\end{proof}

In light of this result, one might speculate that the
$\A$--nonresonance condition \eqref{eqn:stvImproved} and the
$(\A,H)$--nonresonance condition \eqref{eqn:nonres1} are equivalent. 
This is not the case, as the following examples illustrate.  For
simplicity, these examples involve rank $1$ local systems.  In this
context, it is customary to refer to the collection of endomorphisms
$(P_1,\dots,P_m)$ as {\em weights}, and to write
$\la=(\la_1,\dots,\la_m)=(P_1,\dots,P_m)$.

\begin{exm}
Let $\A$ be the arrangement of five lines in $\bP^2$ defined by
the polynomial $Q=x(x-z)y(y-z)z$.  Order the hyperplanes of $\A$
as indicated by the order of the factors of $Q$.  The dense edges
of $\A$ are the hyperplanes, $X_{125}=H_1 \cap H_2 \cap H_5$, and
$X_{345}$. For this arrangement, the weights
$\la=\left(\frac{1}{2},\frac{1}{2},\frac{1}{2},\frac{1}{2},-2\right)$
satisfy the $\A$--nonresonance condition \eqref{eqn:stvImproved},
but there is no integer vector $k$ for which $\la+k$ satisfies the
$(\A,H)$--nonresonance condition \eqref{eqn:nonres1}.
\end{exm}

Note that the rank $1$ local system $\LL$ corresponding to the weights
$\la$ in the previous example has trivial monodromy about one of the
hyperplanes of $\A$.  There are examples where the monodromy about
each hyperplane is nontrivial, and \eqref{eqn:stvImproved} holds, but
\eqref{eqn:nonres1} does not hold for any hyperplane of $\A$.

\begin{exm}
Let $\A$ be the arrangement of six lines in $\bP^2$ defined by the
polynomial $Q=x(x-z)y(y-2z)(x-y)z$.  Order the hyperplanes of $\A$
as indicated by the order of the factors of $Q$.  The dense edges
of $\A$ are the hyperplanes, $X_{126}$, $X_{135}$, and $X_{346}$.
For this arrangement, the weights
$\la=\left(-\frac{5}{3},\frac{1}{3},-\frac{5}{3},
\frac{1}{3},\frac{7}{3},\frac{1}{3}\right)$ satisfy the
$\A$--nonresonance condition \eqref{eqn:stvImproved}, but there is
no integer vector $k$ for which $\la+k$ satisfies the
$(\A,H)$--nonresonance condition \eqref{eqn:nonres1}.
\end{exm}

In both examples, the monodromy about each rank 2 dense edge is
trivial, so the $(\A,H)$--nonresonance condition \eqref{eqn:nonres1}
cannot hold for any integer translate of the weights, but the weights
nevertheless satisfy the $\A$--nonresonance condition
\eqref{eqn:stvImproved}.

\section{An application to Milnor fibers of line arrangements}
\label{sec:mf}

Let $m=|\A|$ be the number of hyperplanes in the arrangement 
$\A \subset \bP^n$ and let $Q=0$ be a reduced equation for the
corresponding central arrangement $\widetilde\A$ in $\C^{n+1}$.  The
smooth affine hypersurface $F$ in $\C^{n+1}$ given by the equation
$Q=1$ is called the Milnor fiber of the arrangement $\A$ (resp.,
$\widetilde\A$).  There is a naturally associated monodromy operator
$h : F \to F$ given by multiplication by $\tau = \exp(2\pi \ii/m)$,
satisfying $h^m=1$.  For $k,p \in \Z$, let $b_p(F)_k$ denote the
dimension of the $\tau^k$-eigenspace corresponding to the monodromy
action on $H^p(F,\C)$.  It is known that
\[
b_p(F)_k=\dim H^p(M(\A),\LL _k)
\]
where $\LL _k$ is the rank one local system on $M(\A )$ corresponding
to the monodromies $T_1=T_2= \dots =T_m=\tau ^k$, see \cite{CS}.  In
this section, we prove the following strengthening of a result of
Massey \cite{Ma}.

\begin{thm} \label{thm:mf}
Let $\A$ be a line arrangement in $\bP ^2$, with associated Milnor
fiber $F$.  Then for any integer $0<k <m$ and any line $H$ in the
arrangement $\A$ we have
\[
b_1(F)_k \leq \sum _x (m_x-2)
\]
where the sum is over all points $x \in H$ such that the multiplicity
of $\A$ at $x$ is $m_x >2$ and $m$ divides $km_x$.
\end{thm}
\begin{proof}
In this proof, we work directly on $\bP^2$ (without using the partial
resolution $W$ from Section \ref{sec:proof}).  Let
$i: M(\A) \to \bP^2 \setminus H$ and $j: \bP^2 \setminus H \to \bP^2$
denote the natural inclusions, and set $\F=Ri_*\LL[2]$.  Recall that
$\F$ is perverse, and extend the canonical morphism
$Rj_!\F \to Rj_*\F$ in $\b{D}_c^b(\bP^2)$ to a distinguished triangle
\begin{equation} \label{eqn:Delta}
Rj_!\F  \to Rj_*\F  \to \GG  \to
\end{equation}

Let $x$ be a point of multiplicity $m_x$ on the chosen line $H$. 
Applying the functors $\HH^k _x$ to the triangle \eqref{eqn:Delta}, we
obtain a long exact sequence, which yields in particular
\[
\HH^{-1}(\GG)_x =\HH^{-1}(Rj_*\F)_x=H^{1}(M_x,\LL)
\text{ and }
\HH^{0}(\GG)_x =\HH^{0}(Rj_*\F)_x=H^{2}(M_x,\LL).
\]
Here $M_x =M(\A) \cap B_x $, where $B_x$ is a small open ball
centered at $x$. It follows that $M_x$ is homeomorphic to the
complement of the central line arrangement in $\C^2$ defined by
$f(y)=y_1^{m_x}+y_2^{m_x}$. There are two cases. If $m$ does not
divide $km_x$, then $H^q(M_x,\LL)=0$ for all $q$. On the other
hand, if $m$ divides $km_x$, then
$\dim H^{1}(M_x,\LL)=\dim H^{2}(M_x,\LL)=m_x-2$, see for instance
\cite[p.~109]{Di} and \cite{CS,Ma}. 

Now apply the functors $\bH^k$ to the triangle \eqref{eqn:Delta},
and use the vanishing results \eqref{eqn:vanishing}. This yields
an exact sequence
\[
0  \to H^{1}(M(\A),\LL)  \to \bH^{-1}(\bP ^2, \GG)  \to
H_c^{2}(M(\A),\LL) \to \dots
\]
To compute the middle term in this sequence, we use the spectral
sequence
\[
E_1^{p,q}=H^p(S,\HH^q \GG)
\]
where $S$ is the support of $\GG$, a finite set.  It follows that
$\bH^{-1}(\bP^2,\GG)=H^0(S,\HH^{-1} \GG)$ is a $\C$-vector space of
dimension $\sum_x (m_x-2)$ where the sum is over all points $x \in H$
such that the multiplicity of $\A$ at $x$ is $m_x >2$ and $m$ divides
$km_x$.
\end{proof}
\begin{rem}
The proof of Theorem \ref{thm:mf} uses only the fact that the support
of the sheaf $\GG$ is finite.  This may happen for arrangements in
$\bP^n$ for $n>2$, yielding a more general version of the theorem.
\end{rem}

\section{Generalizations} \label{sec:gen}

Our main results can be stated (and proved in the same way, either
using partial resolutions or working directly in the projective space)
for arrangements of hypersurfaces $N= \bigcup_{i=1}^{m}V_i$ in a
projective space $\bP^n$.  See \cite{Da} and the references there for
other results in this setting.  It is not necessary to assume that the
individual hypersurfaces are smooth.  It is enough to impose local
vanishing assumptions, both for the intersections contained in a fixed
hypersurface, say $V_1$, and at all singular points of $V_1$ itself.

In the hyperplane arrangement case, we can treat the local cohomology
groups whose vanishing is necessary in the proof of Lemma
\ref{lem:canonical} in terms of complements of central arrangements. 
This allows us to decrease the dimension by one and proceed by
induction.

In the general hypersurface arrangement case, this induction is no
longer available, since the complement need not be locally a cone over
a projective arrangement of smaller dimension.  However, the following
approach may be used to obtain vanishing results in this generality.

Let $f:(\C^{n+1},0) \to (\C,0)$ be an analytic function germ and let
$\F=\C[n+1]$ be the perverse sheaf on $\C^{n+1}$ obtained by shifting
the constant sheaf $\C$.  It is known that perverse sheaves are
preserved by the perverse vanishing cycle functor, 
\cite[Corollary 10.3.13]{KS}.  Thus ${\,}^p \phi_f (\F) \in \Perv(X)$,
where $X=f^{-1}(0)$.  There is a natural monodromy automorphism $\mu:
{\,}^p \phi_f (\F) \to {\,}^p \phi_f (\F)$.  For any $a \in \C$, we
can consider the eigenspace $\F_a = \ker(\mu -a \cdot \hbox{Id})$,
which is a well-defined perverse sheaf on $X$, since the category
$\Perv(X)$ is abelian, \cite[Proposition 10.1.11]{KS}.

For any point $x \in X$, $\HH^m({\,}^p \phi_f (\F))_x=H^{m+n}(F_x)$,
where $F_x$ is the local Milnor fiber of $f$ at the point $x$. 
Moreover, the induced action of $\mu$ on $\HH^m({\,}^p \phi_f (\F))_x$
corresponds exactly to the usual monodromy action on the local Milnor
fiber $F_x$.

Let $S_a$ be the support of the sheaf $\F_a$ and let $s_a=\dim S_a$,
with the convention $\dim \emptyset =-1$.  Note that the integer $s_a$
depends only on the hypersurface germ $(X,0)$: indeed, any two reduced
equations for this germ are topologically equivalent (since the
contact equivalence classes are connected), 
see \cite[Remark 3.1.8]{Di}.

It follows that $\F_a \in \Perv(S_a)$, see \cite[Section (5.2)]{Di2}. 
Hence the support condition in the definition of perverse sheaves
gives $\HH^m({\,}^p \phi_f (\F))_x=0$ for any $m < -s_a$.  This
implies that
\[
H^{n-s_a-j}(F_0)_a=0
\]
for all $j>0$. 

Using the Milnor fibration of $f$ at the origin, we can identify the
corresponding Milnor fiber $F_0$ with an infinite cyclic covering of
$U_0$, the local complement of $X$ in $(\C^{n+1},0)$.  For $a \in
\C^*$, we denote by $\LL_a$ the rank one local system on $U_0$ whose
monodromy around each irreducible component of $X$ is multiplication
by $a$.

Then it is well known that
\[
\dim H^q(U_0,\LL_a)=\dim H^{q-1}(F_0)_a + \dim H^{q}(F_0)_a,
\]
see for instance  \cite{Li2,DN}. It follows that
\[
\dim H^q(U_0,\LL_a)=0,
\]
for all $q \leq n-1-s_a$.  Applying this local vanishing result to the
global setting of hypersurface arrangements as in Section
\ref{sec:proof} above, we obtain the following (note that $n+1$ is
replaced by $n$!).

\begin{thm} 
Let $N= \bigcup_{i=1}^m V_i$ be a hypersurface arrangement in
$\bP^n$, with associated Milnor fiber $F$.  Let $d=d_1+\cdots +d_m$ be
the degree of $N$.  For each point $x \in V_1$, denote by $s(x,k)$ the
number $s_{\tau ^k}$ associated to the hypersurface germ $(N,x)$ as
above, with $\tau =\exp(2\pi \ii/d)$.  Let $s_k=\max_{x \in
V_1}s(x,k)$.  Then for any integer $0<k <d$, we have
\[
b_q(F)_k =0
\]
for all $q \leq n-2-s_k$.
\end{thm}

\begin{cor}
If $N= \bigcup_{i=1}^m V_i$ is a normal crossing divisor at any point
$x \in V_1$, then the monodromy action on $H^q(F)$ is trivial for $q
\leq n-1$.
\end{cor}

\bibliographystyle{amsalpha}

\end{document}